\documentclass [11pt,oneside,a4paper,mathscr]{amsart}

\usepackage{amscd,amsmath,amssymb,euscript}
\usepackage[frame,cmtip,curve,arrow,matrix,line,graph]{xy}

\title{Stability conditions on triangulated categories}
\author{Tom Bridgeland}

\date{}

\newtheorem{thm}{Theorem}[section]

\newtheorem{cor}[thm]{Corollary}
\newtheorem{prop}[thm]{Proposition}

\newtheorem{lemma}[thm]{Lemma}

\newenvironment{pf}{\paragraph{Proof}}{\qed\par\medskip}

\renewcommand{\leq}{\leqslant}
\renewcommand{\geq}{\geqslant}

\theoremstyle{definition}
\newtheorem{defn}[thm]{Definition}

\newtheorem{example}[thm]{Example}

\newcommand{\Slice}{\operatorname{Slice}}

\newcommand{\half}{\frac{1}{2}}

\newcommand{\quarter}{\frac{1}{4}}

\newcommand{\K}{{{K}}}

\newcommand{\Aut}{\operatorname{Aut}}
\newcommand{\U}{\mathcal{U}}
\newcommand{\T}{\operatorname{\mathcal D}}

\newcommand{\isom}{\cong}

\newcommand{\tensor}{\otimes}
\newcommand{\PP}{\operatorname{\mathbb P}}
\newcommand{\M}{\operatorname{\mathcal M}}
\newcommand{\C}{\mathbb C}

\newcommand{\QQ}{\mathbb Q}

\newcommand{\F}{\mathcal F}
\newcommand{\G}{\mathcal G}
\newcommand{\Z}{\mathbb Z}
\newcommand{\A}{\mathcal A}
\newcommand{\B}{\mathcal B}
\newcommand{\CC}{\mathcal C}
\newcommand{\OO}{\mathcal O}

\newcommand{\onto}{\twoheadrightarrow}
\renewcommand{\P}{\mathcal P}
\newcommand{\Q}{\mathcal Q}

\newcommand{\D}{\operatorname{D}}

\newcommand{\Hom}{\operatorname{Hom}}

\newcommand{\SL}{\operatorname{SL}}
\newcommand{\eu}{\operatorname{\chi}}

\newcommand{\coim}{\operatorname{coim}}
\renewcommand{\Re}{\operatorname{Re}}
\renewcommand{\Im}{\operatorname{Im}}
\newcommand{\lRa}[1]{\xrightarrow{\ #1\ }}

\newcommand{\lra}{\longrightarrow}

\newcommand{\rk}{\operatorname{rank}}

\newcommand{\Stab}{\operatorname{Stab}}
\newcommand{\R}{\mathbb{R}}

\newcommand{\GL}{\operatorname{GL^+}}

\newcommand{\ha}{\frac{1}{2}}
\newcommand{\N}{\operatorname{\mathcal{N}}}

\newcommand{\ZZ}{\operatorname{\mathcal{Z}}}
\newcommand{\coker}{\operatorname{coker}}
\newcommand{\im}{\operatorname{im}}
\newcommand{\grp}{{\tilde{\GL}}(2,\R)}

\newcommand{\m}{m}

\begin{document}

\maketitle

\section{Introduction}
This paper introduces the notion of a stability condition on a
triangulated category.
The motivation comes from
the study of Dirichlet branes in string theory,
and especially from M.R. Douglas's work on $\Pi$-stability.
From a mathematical point of view,
the most interesting feature of the definition is that
the set of stability conditions $\Stab(\T)$ on a fixed category $\T$
has a natural topology, thus defining a new invariant of triangulated
categories.
In a separate article I shall give a detailed description of this space
of stability conditions when $\T$ is the bounded derived category of
coherent sheaves on a K3 surface \cite{Br}. The present paper though is
almost pure homological algebra. After setting up the necessary
definitions I prove a deformation result which
shows that the space $\Stab(\T)$ with its natural topology
is a manifold, possibly  infinite-dimensional.
 
\subsection{}
Before going any further let me
describe a simple example of a stability condition on a
triangulated category.
Let $X$ be a nonsingular projective curve
and let $\D(X)$ be its bounded derived category of
coherent sheaves.
Recall \cite{HN} that any nonzero coherent sheaf $E$ on $X$
has a unique Harder-Narasimhan filtration
\[ 0=E_0 \subset E_1\subset \cdots \subset E_{n-1}\subset E_n=E, \]
whose factors $E_j/E_{j-1}$ are semistable sheaves with descending
slope $\mu=\deg/\rk$. Torsion sheaves should be thought of as having slope
$+\infty$ and come first in the filtration.
On the other hand, given an object $E\in\D(X)$,
the truncations $\sigma_{\leq j}(E)$ associated
to the standard t-structure on 
$\D(X)$ fit into triangles
\[\xymatrix@C=.4em{
\ldots \ar[rr] && \sigma_{\leq j-1}(E) \ar[rrrr] &&&& \sigma_{\leq
j}(E) \ar[rrrr] \ar[dll] &&&& \sigma_{\leq j+1}(E)
\ar[rr] \ar[dll] && \ldots \\
&&&& A_j \ar@{-->}[ull] &&&& A_{j+1} \ar@{-->}[ull]
}\]
which allow one to break up $E$ into
its shifted cohomology sheaves $A_j=H^j(E) [-j]$.
Combining these two ideas, one can concatenate the Harder-Narasimhan
filtrations of the cohomology sheaves $H^j(E)$ to obtain a
kind of filtration of any nonzero object $E\in\D(X)$ by
shifts of
semistable sheaves.

Now define a complex-valued linear map on the Grothendieck group
$\K(X)$ of $\D(X)$ by the formula
\[Z(E)=-\deg(E)+i\rk(E).\]
For each nonzero sheaf $E$ on $X$, there is a
unique branch $\phi(E)$ of
$(1/\pi)\arg Z(E)$
lying in the interval $(0,1]$. If one defines
\[\phi\big(E[k]\big)=\phi(E)+k,\]
for each integer $k$, then
the filtration described above is by objects of descending
phase $\phi$, and in fact is unique with this property.
Thus each nonzero object of $\D(X)$ has a kind of generalised
Harder-Narasimhan filtration.  Note that not all objects of $\D(X)$ have
a well-defined phase, indeed many objects of $\D(X)$ define the zero
class in $\K(X)$. Nonetheless, the phase function is well-defined on
the generating subcategory
$\P\subset\D(X)$ consisting of shifts of semistable sheaves, 
and in fact defines an $\R$-grading of this category.  

\subsection{}
The definition of a stability condition on a triangulated category is
obtained by abstracting these generalised Harder-Narasimhan
filtrations
of nonzero objects of $\D(X)$ together with
the map $Z$ as follows. Throughout the paper the Grothendieck group of a triangulated category $\T$ is denoted $\K(\T)$.

\begin{defn}
\label{pemb}
A stability condition $(Z,\P)$ on a triangulated category $\T$
consists of
a group homomorphism
$Z\colon\K(\T)\to\C$ called the \emph{central charge},
and full additive
subcategories $\P(\phi)\subset\T$ for each $\phi\in\R$,
satisfying the following axioms:
\begin{itemize}
\item[(a)] if $E\in \P(\phi)$ then $Z(E)=\m(E)\exp(i\pi\phi)$ for some
 $\m(E)\in\R_{>0}$,
\item[(b)] for all $\phi\in\R$, $\P(\phi+1)=\P(\phi)[1]$,
\item[(c)] if $\phi_1>\phi_2$ and $A_j\in\P(\phi_j)$ then $\Hom_{\T}(A_1,A_2)=0$,
\item[(d)] for each nonzero object $E\in\T$ there is a finite sequence of real
numbers
\[\phi_1>\phi_2> \cdots >\phi_n\]
and a collection of triangles
\[
\xymatrix@C=.5em{
0_{\ } \ar@{=}[r] & E_0 \ar[rrrr] &&&& E_1 \ar[rrrr] \ar[dll] &&&& E_2
\ar[rr] \ar[dll] && \ldots \ar[rr] && E_{n-1}
\ar[rrrr] &&&& E_n \ar[dll] \ar@{=}[r] &  E_{\ } \\
&&& A_1 \ar@{-->}[ull] &&&& A_2 \ar@{-->}[ull] &&&&&&&& A_n \ar@{-->}[ull] 
}
\]
with $A_j\in\P(\phi_j)$ for all $j$.

\end{itemize}
\end{defn}

I shall always assume that the category $\T$ is essentially small,
that is, that $\T$ is equivalent to a category in which the
class of objects is a set. One can then consider the set
of all stability conditions on $\T$. In fact it makes more sense to
restrict attention to stability conditions satisfying
a certain technical
condition called local-finiteness (Definition \ref{lofo}).
I show how to put a natural
topology on the set $\Stab(\T)$ of such stability conditions
and prove the following theorem.

\begin{thm}
\label{lasty}
Let $\T$ be a triangulated category.
For each connected component $\Sigma\subset\Stab(\T)$
there is a linear subspace $V(\Sigma)\subset \Hom_{\Z}(\K(\T),\C)$
with a well-defined linear topology and a local homeomorphism
$\ZZ\colon\Sigma\to V(\Sigma)$ which maps a stability condition $(Z,\P)$
to its central charge $Z$.
\end{thm}

It follows immediately from this theorem that each component
$\Sigma\subset\Stab(\T)$ is a
manifold, locally modelled on the topological vector space $V(\Sigma)$.

\subsection{}
Suppose now that $\T$ is linear over a field $k$. This means
that the morphisms of $\T$ have the structure of a vector space over
$k$, with respect to which the composition law is bilinear. Suppose
further that $\T$ is of finite type, that is
that for every pair of objects
$E$ and $F$ of $\T$ the vector space
$\bigoplus_i \Hom_{\T}(E,F[i])$
is finite-dimensional. In this situation one can define a bilinear form on
$\K(\T)$,
known as the Euler form,
via the formula
\[\eu(E,F)=\sum_i (-1)^i \dim_k \Hom_{\T}(E,F[i]),\]
and a free abelian group
$\N(\T)=\K(\T)/\K(\T)^{\perp}$ called the \emph{numerical Grothendieck group}
of $\T$. If this group $\N(\T)$ has finite rank the category
$\T$ is said to be \emph{numerically finite}.

Suppose then that $\T$ is of finite type over a field, and
define $\Stab_{\N}(X)$ to be
the subspace of $\Stab(\T)$ consisting of
\emph{numerical} stability conditions, that is, those for which the central
charge $Z\colon\K(\T)\to\C$ factors through the quotient group
$\N(\T)$.
The following result is an immediate consequence of Theorem \ref{lasty}.
\begin{cor}
Suppose $\T$ is numerically finite.  For
each connected component $\Sigma\subset \Stab_{\N}(\T)$ there is a
subspace $V(\Sigma)\subset \Hom_{\Z}(\N(\D),\C)$ and a local
homeomorphism
$\ZZ\colon\Sigma\to V(\Sigma)$ which maps a stability condition to its
central charge $Z$. In particular $\Sigma$ is a finite-dimensional
complex manifold.
\end{cor}

There are two large classes of examples of numerically finite
triangulated categories.
Firstly, if $A$ is a finite-dimensional algebra over a field, then the
bounded derived category $\D(A)$ of finite-dimensional left $A$-modules
is numerically finite. The corresponding space of numerical stability
conditions will be denoted $\Stab(A)$. Secondly, if $X$ is a
smooth projective variety over $\C$
then the Riemann-Roch theorem shows that
the bounded derived
category $\D(X)$ of coherent sheaves on $X$
is numerically finite. In this case the space of numerical stability
conditions will be denoted $\Stab(X)$.

Obviously one would like to be able to compute these spaces
of stability conditions in some interesting examples.
The only case considered in this paper
is when $X$ 
is an elliptic curve. Here the
answer is rather straightforward: $\Stab(X)$
is connected, and
there is a local homeomorphism
\[\ZZ\colon\Stab(X)\to \C^2.\]
The image of this map is $\GL(2,\R)$, the group of rank two matrices
with positive determinant, considered as an open subset of $\C^2$ in
the obvious way, and $\Stab(X)$ is the universal cover of this
space. Perhaps of more interest is the quotient of $\Stab(X)$ by the
group of autoequivalences of $\D(X)$. One has
\[{\Stab(X)}\left/{\Aut \D(X)}\right.\isom
{\GL(2,\R)}\left/ \SL(2,\Z),\right.\]
which is a $\C^*$-bundle over the modular curve.

\subsection{}
The motivation for the definition of a stability condition given above
came from the work of Douglas on $\Pi$-stability for Dirichlet
branes. It therefore seems appropriate to include here
a short summary of some of Douglas' ideas. However the author is
hardly an expert in this area, and this section will inevitably contain
various inaccuracies and over-simplifications.
The reader would be well-advised to
consult the original papers of Douglas \cite{Do1,Do2,Do4} and
Aspinwall-Douglas \cite{Do3}.
Of course, those with no interest
in string theory can
happily skip to the next section.

String theorists believe that the supersymmetric
nonlinear sigma model allows them to
associate a $(2,2)$
superconformal field theory (SCFT) to a set of data consisting of
a compact, complex
manifold $X$
with trivial canonical bundle,
a K{\"a}hler class $\omega\in H^2(X,\R)$ and a class
$B\in H^2(X,\R/\Z)$ induced by a closed 2-form
on $X$ known as the B-field. Assume for simplicity
that $X$ is a simply-connected threefold.
The set of possible choices
of this data up to equivalence then defines an open subset
$\U_X$ of the moduli space $\M$ of SCFTs. This moduli space $\M$ has two
foliations, which when restricted to $\U_X$ just correspond to
those obtained by holding constant either the complex structure of $X$ or the
complexified K{\"a}hler class $B+i\omega$.

It is worth bearing in mind that the open subset $\U_X\subset\M$ described
above is just a neighbourhood of a particular `large volume limit' of
$\M$; a given
component of $\M$ may contain points corresponding to sigma
models on topologically
distinct manifolds $X$ and also points which do not correspond
to sigma models at all. One of the long-term goals of the
present work is to try to gain a clearer mathematical
understanding of this moduli space $\M$.

The next step is to consider branes.
These are boundary conditions
in the SCFT and naturally form the objects of a category, with
the space of morphisms
between a pair of branes being the spectrum of open strings with
boundaries on
them.
One of the most striking claims of recent work in string theory
is that the SCFT corresponding to a nonlinear sigma model
admits a `topological twisting' in which the corresponding category of
branes is equivalent to $\D(X)$,
the bounded derived category of coherent sheaves on
$X$. In particular this category does not depend on
the so-called stringy K{\"a}hler moduli space of $X$, that is the leaf
$\M_K(X)\subset\M$
corresponding to a fixed complex structure on $X$.

Douglas starts from this point of view and proceeds to argue that at
each point in $\M_K(X)$ there is a subcategory $\P\subset\D(
X)$ whose objects are the physical or BPS branes for the
corresponding SCFT. He also gives a precise criterion
`$\Pi$-stability' for describing how this subcategory $\P$ changes along
continuous paths in $\M_K(X)$. An important  point to note is that
whilst the category of BPS branes is
well-defined at any point in $\M_K(X)$, the embedding
$\P\subset\D( X)$ is not, so that monodromy around
loops in the K{\"a}hler moduli space leads
to different subcategories $\P\subset\D(X)$,
related to each other
by autoequivalences of $\D(X)$.

The definition of a stability condition given above was an attempt to
abstract the properties of the subcategories $\P\subset\D(X)$. Thus
the points of the K{\"a}hler moduli space $\M_K(X)$ should be thought
of as defining points in the quotient $\Stab(X)/\Aut \D(X)$, and the
category $\P=\bigcup_{\phi}\P(\phi)$ should be thought of
as the category of BPS branes at the corresponding point of
$\M_K(X)$.

There is also a mirror side to this story. According to
the predictions of mirror symmetry there is an involution $\sigma$ of
the moduli space $\M$ which identifies some part of the open subset
$\U_X$ defined
above with part of the corresponding set $\U_{\check{X}}$
associated to a mirror
manifold $\check{X}$. This identification exchanges the two foliations, so that
the K{\"a}hler moduli space of $X$ becomes identified with the moduli
of complex structures on $\check{X}$ and vice versa.

Kontsevich's
homological mirror conjecture \cite{Kon} predicts that
 the derived category $\D(X)$
is equivalent to the derived Fukaya category
$\D\operatorname{Fuk}(\check{X})$.
Roughly speaking, this equivalence is expected to take
the subcategory $\P(\phi)\subset\D(X)$ at a particular point of
$\M_K(X)$ to the subcategory of $\D\operatorname{Fuk}(\check{X})$
consisting of special Lagrangians of phase $\phi$ with respect to the
corresponding complex structure on $\check{X}$.
%However it seems to be completely unknown whether objects of the
%derived Fukaya category have Harder-Narasimhan-type filtrations by
%special Lagrangians of descending phase as in Definition \ref{pemb}.
For more on this side of the picture see for example \cite{Th1,Th2}.

\subsection*{Notation}
The term generalised metric will be used to mean a distance function
$d\colon X\times X\to [0,\infty]$ on a set $X$ satisfying all the usual
metric space axioms except that it need not be finite.
Any such function defines a topology on $X$ in the usual way
and induces a
metric space structure
on each connected component of $X$.

The reader is referred to \cite{GM,Ht,Ve} for
background on triangulated categories. I always assume that my
categories are essentially small.
I write [1] for the shift (or
translation) functor of a triangulated category and draw my triangles as
follows
\[\xymatrix@C=.5em{A \ar[rrrr] &&&& B  \ar[dll] \\
 && C \ar@{-->}[ull]}\]
where the dotted arrow means a morphism $C\to A[1]$. Sometimes I just
write \[A\lra B\lra C.\]

The Grothendieck group of a triangulated category $\T$ is denoted $\K(\T)$. Similarly, the Grothendieck group of an abelian category $\A$ is denoted $\K(\A)$.

A full subcategory $\A$ of a triangulated category $\T$
will be called extension-closed if
whenever $A\to B\to C$
 is a triangle in $\T$ as above, with $A\in\A$ and $C\in\A$,
then $B\in\A$ also.
The extension-closed subcategory of $\T$ generated by a
full subcategory $\mathcal{S}\subset\T$ is 
the smallest extension-closed
full subcategory of $\T$ containing $\mathcal{S}$.

\subsection*{Acknowledgements} My main debt is to Michael Douglas
whose papers on $\Pi$-stability provided the key idea for this
paper. I'm also indebted to Dmitry Arinkin and Vladimir Drinfeld
who pointed out a simpler way to prove Theorem \ref{biggy}.
Finally I'd like to thank
Alexei Bondal, Mark Gross, Alastair King,
Antony Maciocia, So Okada,
Aidan Schofield and Richard Thomas for their comments and corrections.
% ********************************************************************
% ********************************************************************
% ********************************************************************

\section{Stability functions and Harder-Narasimhan filtrations}

The definition of a stable vector bundle on a curve
has two fundamental ingredients, namely the partial ordering
$E\subset F$
arising from the notion of a sub-bundle, and the numerical
ordering coming from the slope function $\mu(E)$.
Both of these ingredients were
generalised by A.N. Rudakov \cite{Ru2} to give an abstract notion of a
stability condition on an abelian category.
For the purposes of this paper, it will not be necessary to
adopt the full generality of Rudakov's
approach, which allowed for arbitrary orderings on abelian categories.
In fact one need only  consider
orderings induced by certain phase functions, as follows.

\begin{defn}
A \emph{ stability function} on
an abelian  category $\A$ is a group homomorphism $Z\colon \K(\A)\to\C$
such that for all $0\neq E\in\A$ the complex number $Z(E)$ lies in the
strict upper half-plane
$H=\{r\exp(i\pi\phi):r>0\text{ and }0<\phi\leq 1\}
\subset\C.$
\end{defn}

Given a stability function $Z\colon\K(\A)\to\C$,
the \emph{phase} of an object $0\neq E\in\A$ is defined to be
\[\phi(E)=(1/\pi)\arg Z(E)\in (0,1].\]
The function $\phi$ allows one to order the nonzero objects of
the category $\A$ and thus leads to a notion of stability for objects
of $\A$.
Of course one could equally well define this ordering
using the function $-\Im Z(E)/\Re
Z(E)$ taking values in $(-\infty,+\infty]$, but
in what follows it will be important to use
the phase function $\phi$
instead.

\begin{defn}
Let $Z\colon\K(\A)\to\C$ be a  stability function on an abelian
category $\A$.
An object $0\neq
E\in\A$ is said to be \emph{semistable} (with respect to $Z$)
if every subobject $0\neq A\subset E$ satisfies $\phi(A)\leq\phi(E)$.
\end{defn}

Of course one could equivalently define a semistable object
$0\neq E \in\A$ to be one for which $\phi(E)\leq\phi(B)$
for every nonzero quotient $E\onto B$. The
importance of semistable objects in this paper
is that they provide a way to filter
objects of $\A$. This is the so-called Harder-Narasimhan property,
which was first proved for bundles on curves in \cite{HN}.

\begin{defn}
\label{HN}
Let $Z\colon\K(\A)\to\C$ be a  stability function on an abelian
category $\A$. A \emph{Harder-Narasimhan filtration} of an
object $0\neq E\in\A$ is a
finite chain of subobjects
\[0=E_0\subset E_1\subset \cdots\subset E_{n-1}\subset E_n=E\]
whose factors $F_j=E_j/E_{j-1}$ are semistable objects of $\A$ with
\[\phi(F_1)>\phi(F_2)>\cdots>\phi(F_n).\]
The stability function $Z$ is said to have the \emph{Harder-Narasimhan property} if
every nonzero object of $\A$ has a Harder-Narasimhan filtration.
\end{defn}

Note that
if $f\colon E\to F$ is a nonzero map between
semistable objects
then by considering $\im f\isom\coim f$ in the usual way, one sees
that $\phi(E)\leq\phi(F)$. It follows easily from this that
Harder-Narasimhan
filtrations (when they exist) are unique.
The following slight strengthening of a result of
Rudakov \cite{Ru2} shows that the existence of Harder-Narasimhan
filtrations is actually a
rather weak assumption.

\begin{prop}
\label{rud}
Suppose $\A$ is an abelian category with a 
stability function
$Z\colon \K(\A)\to \C$ satisfying the chain conditions
\begin{itemize}
\item[$(a)$] there are no infinite sequences
of subobjects in $\A$
\[\cdots\subset E_{j+1}\subset E_{j}\subset\cdots\subset E_2\subset
E_1\]
with $\phi(E_{j+1})>\phi(E_j)$ for all $j$,

\item[$(b)$] there are no infinite sequences
of quotients in $\A$
\[E_1\onto E_2\onto \cdots\onto E_j\onto E_{j+1}\onto\cdots\]
with $\phi(E_{j})>\phi(E_{j+1})$ for all $j$.
\end{itemize}
Then $\A$ has the Harder-Narasimhan property.
\end{prop}

\begin{pf}
First note that if $E\in\A$ is nonzero then  either $E$ is
semistable or there is a subobject $0\neq E'\subset E$ with
$\phi(E')>\phi(E)$. Repeating the argument and using the first chain
condition it follows that every nonzero object of $\A$ has a semistable
subobject $A\subset E$ with $\phi(A)\geq\phi(E)$.
A similar argument using the second chain condition gives the
dual statement: every nonzero object of $\A$
has a semistable quotient $E\onto B$ with $\phi(E)\geq\phi(B)$.

A maximally destabilising quotient (mdq) of an object $0\neq E\in\A$
is defined to be
a nonzero quotient $E\onto B$ such that any
nonzero quotient $E\onto B'$ satisfies $\phi(B')\geq\phi(B)$, with
equality holding only if $E\onto B'$ factors via $E\onto B$.
By what was said
above it is enough to check this condition under the additional
assumption that $B'$ is semistable.
Note also that if $E\onto B$ is a mdq then $B$ must be semistable with
$\phi(E)\geq\phi(B)$. The first step in the proof of the Proposition
is to show that mdqs always exist.

Take a nonzero object $E\in\A$.
Clearly if $E$ is semistable then the identity map $E\to E$ is a
mdq. Otherwise, as above, there is a short exact sequence
\[0\lra A\lra E\lra E'\lra 0\]
with $A$ semistable and $\phi(A)>\phi(E)>\phi(E')$. I claim that if
$E'\onto B$ is a mdq
for $E'$ then the induced quotient $E\onto B$ is a mdq for
$E$. Indeed, if $E\onto B'$ is a quotient with $B'$ semistable and
$\phi(B')\leq\phi(B)$ then
$\phi(B')<\phi(A)$ so there is no map $A\to B'$ and the quotient
$E\onto B'$ factors via $E'$, which proves the claim.
Thus I can replace $E$ by $E'$ and repeat the argument.
By the second chain condition, this process must eventually  terminate.
It follows that
every nonzero object of $\A$ has a mdq.

Take a nonzero object $E\in\A$. If $E$ is semistable then $0\subset E$
is a Harder-Narasimhan filtration of $E$.
Otherwise there is a short exact sequence
\[0\lra E'\lra E\lra B\lra 0\]
with $E\onto B$ a mdq and $\phi(E')>\phi(E)$.
Suppose $E'\onto B'$ is a mdq. Consider the
following diagram of short exact sequences
\begin{equation}
\begin{CD}
@. @. 0 @. 0 @.\\
@.@. @VVV @VVV\\
0 @>>> K @>>> E' @>>> B' @>>> 0 \\
@. @| @VVV @VVV @. \\
0 @>>> K @>>> E @>>> Q @>>> 0  \\
@.@. @VVV @VVV @.\\
@.@. B @= B @.\\
@.@.@VVV @VVV@.\\
@. @. 0 @. 0 @.
\end{CD}\tag{$\dagger$}\end{equation}
It follows from the definition of $B$ that $\phi(Q)>\phi(B)$ and hence
$\phi(B')> \phi(B)$.
Replacing $E$ by $E'$ and repeating the process,
one obtains a sequence of subobjects of $E$
\[\cdots \subset E^i\subset E^{i-1}\subset\cdots\subset E^1\subset
E^0=E\]
such that $\phi(E^i)>\phi(E^{i-1})$ and with semistable
factors $F^i=E^i/E^{i-1}$ of ascending phase.
This sequence must terminate by
the first chain condition, and renumbering gives a Harder-Narasimhan
filtration of $E$.
\end{pf} 

% ********************************************************************
% ********************************************************************
% ********************************************************************

\section{T-structures and slicings}

The notion of a
t-structure
was introduced by A. Beilinson,
J. Bernstein and P. Deligne \cite{BBD}. T-structures are the tool which allows
one to see the different abelian categories embedded in a given
triangulated category.
A slightly different
 way to think about t-structures is that they provide a way to
break up objects of a triangulated category into pieces (cohomology objects)
indexed by the integers. The aim of this section is to introduce the notion of a
slicing, which alows one
to break up objects of the category
into finer pieces indexed by the real
numbers.
I start by recalling the definition of a t-structure.

\begin{defn}
A \emph{t-structure} on a triangulated category $\T$
is a full subcategory
$\F\subset \T$, satisfying  $\F[1]\subset\F$,
such that if one
defines
\[\F^{\perp}
=\{G\in \T:\Hom_{\T}(F,G)=0 \text{ for
all }F\in\F\},\]
then
for every
object $E\in\T$ there is a triangle $F\to E\to G$ in $\T$
%\[\xymatrix@C=.5em{F \ar[rrrr] &&&& E  \ar[dll] \\
% && G \ar@{-->}[ull]} \]
with $F\in\F$ and $G\in \F^{\perp}$.
\end{defn}

The motivating example is the \emph{standard t-structure} on the
bounded derived category $\D(\A)$ of an abelian category $\A$,
obtained by taking $\F$ to consist of all those objects of $\D(\A)$
whose cohomology objects $H^i(E)\in\A$ are zero for all
$i>0$.

The \emph{heart} of a t-structure $\F\subset\T$ is the full
subcategory
\[\A=\F\cap\F^{\perp}[1]\subset \T.\]
It was proved in \cite{BBD} that $\A$ is an abelian category, with the
short exact sequences in $\A$ being precisely the
triangles in $\T$ all of whose vertices are objects of $\A$.

A t-structure $\F\subset\T$ is said to be \emph{bounded} if
\[\T=\bigcup_{i,j\in\Z}\F[i]\cap\F^{\perp}[j].\]
A bounded t-structure $\F\subset\T$
is determined by its heart $\A\subset\T$. In
fact $\F$ is the extension-closed subcategory generated by the
subcategories $\A[j]$ for integers $j\geq 0$.
The following easy result gives another characterisation of
bounded t-structures. The proof is a good exercise in manipulating
the definitions.

\begin{lemma}
\label{ll}
Let $\A\subset\T$ be a full additive
subcategory of a triangulated category
$\T$.
Then $\A$ is the heart of a bounded t-structure $\F\subset\T$ if and
only if the
following two conditions hold:
\begin{itemize}
\item[(a)]if $k_1>k_2$ are integers and $A$ and $B$ are objects of
$\A$ then $\Hom_{\T}(A[k_1],B[k_2])=0$,

\item[(b)]
for every nonzero object $E\in\T$ there is a finite
sequence of integers
\[k_1>k_2>\cdots>k_n\]
and a collection of triangles
\[
\xymatrix@C=.5em{
0_{\ } \ar@{=}[r] & E_0 \ar[rrrr] &&&& E_1 \ar[rrrr] \ar[dll] &&&& E_2
\ar[rr] \ar[dll] && \ldots \ar[rr] && E_{n-1}
\ar[rrrr] &&&& E_n \ar[dll] \ar@{=}[r] &  E_{\ } \\
&&& A_1 \ar@{-->}[ull] &&&& A_2 \ar@{-->}[ull] &&&&&&&& A_n \ar@{-->}[ull] 
}
\]
with $A_j\in\A[k_j]$ for all $j$.\qed
\end{itemize}
\end{lemma}

Taking Lemma \ref{ll} as a guide, one can now replace
the integers $k_j$ with real numbers $\phi_j$
to give the notion of a
slicing. This is the key ingredient in the definition of a
stability condition on a triangulated category.
Some explicit examples will be given in Section \ref{stabilityconditions}.

\begin{defn}
\label{slicing}
A \emph{slicing} $\P$ of
a triangulated category $\T$ consists of
full additive subcategories
$\P(\phi)\subset\T$ for each $\phi\in\R$ satisfying the following axioms:

\begin{itemize}
\item[(a)] for all $\phi\in\R$, $\P(\phi+1)=\P(\phi)[1]$,
\item[(b)] if $\phi_1>\phi_2$ and $A_j\in\P(\phi_j)$ then $\Hom_{\T}(A_1,A_2)=0$,
\item[(c)] for each nonzero object $E\in\T$ there is a finite sequence of real
numbers
\[\phi_1>\phi_2> \cdots >\phi_n\]
and a collection of triangles
\[
\xymatrix@C=.5em{
0_{\ } \ar@{=}[r] & E_0 \ar[rrrr] &&&& E_1 \ar[rrrr] \ar[dll] &&&& E_2
\ar[rr] \ar[dll] && \ldots \ar[rr] && E_{n-1}
\ar[rrrr] &&&& E_n \ar[dll] \ar@{=}[r] &  E_{\ } \\
&&& A_1 \ar@{-->}[ull] &&&& A_2 \ar@{-->}[ull] &&&&&&&& A_n \ar@{-->}[ull] 
}
\]
with $A_j\in\P(\phi_j)$ for all $j$.
\end{itemize}
\end{defn}

Let $\P$ be a slicing of a triangulated category $\T$.
It is an easy exercise to check that the decompositions of
axiom (c) are uniquely defined up to isomorphism.
Given a nonzero object $0\neq E\in\T$ define real numbers
$\phi^+_{\P}(E)=\phi_1$ and
$\phi^-_{\P}(E)=\phi_n$. One has an inequality
$\phi^-_{\P}(E)\leq\phi^+_{\P}(E)$ with equality
holding precisely when $E\in\P(\phi)$ for some $\phi\in\R$.
When the slicing $\P$ is clear from the context I often drop it from the
notation and write $\phi^{\pm}(E)$.

For any interval $I\subset\R$, define $\P(I)$ to be the extension-closed
subcategory of $\T$ generated by the subcategories $\P(\phi)$ for
$\phi\in I$. Thus, for example, the full subcategory $\P((a,b))$
consists of the zero objects of $\T$ together with those
objects $0\neq E\in\T$ which satisfy
$a<\phi^-(E)\leq\phi^+(E)<b$.

\begin{lemma}
\label{blgl}
Let $\P$ be a slicing of a triangulated category $\T$ and
let $I\subset \R$ be an interval of length at most one.
Suppose
\[\xymatrix@C=.5em{A \ar[rrrr] &&&& E  \ar[dll] \\
 && B \ar@{-->}[ull]} \]
is a triangle in $\T$, all of whose vertices
are nonzero objects of $\P(I)$. 
Then there are inequalities $\phi^+(A)\leq\phi^+(E)$
and $\phi^-(E)\leq\phi^-(B)$.
\end{lemma}

\begin{pf}
It is enough to prove the first inequality, the second then follows in the same way.
One can also assume that $I=[t,t+1]$ for some $t\in\R$. By definition, if
$\phi=\phi^+(A)$
there is an object
$A^+\in\P(\phi)$ with a nonzero morphism $f\colon A^+\to A$.
Suppose for a contradiction that
$\phi>\phi^+(E)$. Then
there are no nonzero morphisms
$A^+\to E$ and so $f$ factors via $B[-1]$. But $B[-1]\in\P(\leq t)$
so this implies that $\phi\leq t$. Since $\phi^+(E)\geq t$ this gives the required
contradiction.
\end{pf}

Let $\P$ be a slicing of a triangulated category $\T$ as above.
For any $\phi\in \R$ one has pairs of orthogonal
subcategories $(\P(\smash{>}{\phi}),\P(\leq \phi))$ and
$(\P(\geq \phi),\P(\smash{<}{\phi}))$.
Note that the subcategories $\P({>}{\phi})$ and
$\P({\geq}{\phi})$ are closed under left
shifts and thus define t-structures\footnote{There is an unavoidable
clash of notation here:
in the standard notation for t-structures $\Hom_{\T}(E,F)$ vanishes
providing
$E\in\T^{\leq k}$ and $F\in\T^{> k}$, but in the notation for stability
$\Hom_A(E,F)$ vanishes for $E$ and $F$ semistable providing $E$ has slope $>k$
and $F$ has slope $\leq k$.}
on $\T$. So for each $\phi\in\R$ there are
t-structures $\P(>\phi)\subset\P(\geq \phi)$ on $\T$, indexed by
the real numbers, which are compatible in the
sense that
\[\phi\geq \psi\implies \P(>\phi)\subset\P(>\psi)\text{ and }
\P(\geq \phi)\subset \P(\geq \psi).\]
Of course one could axiomatise these compatible t-structures to give a
slightly weaker notion than that of a slicing.
Note that the heart of the t-structure $\P(>\phi)$ is the subcategory
$\P((\phi,\phi+1])\subset\T$, and similarly,
the t-structure $\P(\geq \phi)$ has heart
$\P([\phi,\phi+1))$.
As a matter of convention, the heart of the slicing $\P$
is defined to be the abelian subcategory 
$\P((0,1])\subset\T$.

% ********************************************************************
% ********************************************************************
% ********************************************************************

\section{Quasi-abelian categories}

Let $\P$ be a slicing of a triangulated category $\T$. It was
observed in
the last section that for any
real number $\phi$ the full subcategories $\P((\phi,\phi+1])$ and $\P([\phi,\phi+1))$ of
$\T$ are the hearts of t-structures on $\T$ and hence are
abelian. Suppose instead that $I\subset\R$ is an interval of length
$<1$ and consider the corresponding full subcategory $\A=\P(I)\subset\T$.
In general this category $\A$ will not be abelian, but it does have
a natural exact structure \cite{Qu},
obtained by defining a short exact sequence in $\A$ to be
a triangle in $\T$ all of
whose vertices are objects of $\A$.
In fact this exact structure
is intrinsic to $\A$ and can be derived from the fact that
$\A$ is a so-called
\emph{quasi-abelian category}. Although this notion is not strictly necessary
for the proof of Theorem \ref{lasty}, it seems worthwhile to summarise
the basic
definitions concerning quasi-abelian categories, since they
undoubtedly provide the right context for discussing these
subcategories
$\P(I)\subset\T$. At a first reading it might be a good idea to skip
this section, since it is only really used in Section 7.
The main reference for quasi-abelian categories
is J.-P. Schneiders' paper \cite{Sch}, see also
\cite[Appendix B]{VdB}.

Suppose then that $\A$ is an additive category with kernels and cokernels. Note
that any such category has pushouts and pullbacks. Given
a morphism $f\colon E\to F$ in $\A$, the image of $f$ is
the kernel of the canonical map $F\to \coker f$, and the coimage
of $f$ is the cokernel of the canonical map $\ker f\to E$. There
is a canonical map $\coim f\to\im f$, and $f$ is called \emph{strict}
if this map is an isomorphism. An abelian category is by definition
an additive
category with kernels and cokernels in which all morphisms are
strict. The following definition gives a weaker notion.

\begin{defn}
A \emph{quasi-abelian category} is an additive category $\A$ with
kernels and cokernels such that every pullback of a strict
epimorphism is a strict epimorphism, and every pushout of a
strict monomorphism is a strict monomorphism.
\end{defn}

A strict short exact sequence in a quasi-abelian category $\A$ is a
diagram
\begin{equation*}\tag{$*$}0\lra A\lRa{i} B\lRa{j} C\lra 0\end{equation*}
in which $i$ is the kernel of $j$ and $j$ is the cokernel of $i$. It
follows that $i$ is a strict monomorphism and $j$ is a strict
epimorphism. Conversely, if $i\colon A\to B$ is a strict monomorphism,
the cokernel of $i$ is a strict epimorphism $j\colon B\to C$
whose kernel is $i$.
Similarly, a strict epimorphism $j\colon B\to C$ has a kernel $i$
fitting into
a strict short exact sequence as above.
The class of strict monomorphisms (respectively
epimorphisms) is closed under composition, and if
\[\xymatrix@C=3em{A \ar[r]^f\ar[dr]_h & B\ar[d]^g  \\
 &  C } \]
is a commutative diagram, then $h$ a strict monomorphism implies that
$f$ is a strict monomorphism, and similarly, $h$ a strict epimorphism
implies that $g$ is a strict epimorphism.
These facts are enough to show
that a quasi-abelian category together with its class
of strict short exact sequences is an exact category \cite{Qu}.
The Grothendieck group of $\A$
is defined to be the abelian group $\K(\A)$ generated
by the objects of $\A$,
with a relation $[B]=[A]+[C]$ for each strict short exact
sequence $(*)$.

The following characterization of quasi-abelian categories was proved
by Schneiders \cite[Lemma 1.2.34]{Sch}.

\begin{lemma}
\label{schn}
An additive category $\A$ is quasi-abelian if and only if there are
abelian categories $\A^{\sharp}$ and $\A^{\flat}$ and fully faithful
embeddings $\A\subset\A^{\sharp}$ and $\A\subset\A^{\flat}$ such that
\begin{itemize}
\item[(a)] if $A\to E$ is a monomorphism in $\A^{\sharp}$ with
$E\in\A$ then $A\in\A$ also,
\item[(b)] if $E\to B$ is an epimorphism in $\A^{\flat}$ with $E\in\A$
then $B\in\A$ also.
\end{itemize}
If these conditions hold, the strict short exact sequences in $\A$
are precisely those sequences $(*)$ which are exact in both $\A^{\sharp}$
and  $\A^{\flat}$.
\qed\end{lemma}

A good example to bear in mind is the category $\A$ of torsion-free
sheaves on a smooth projective variety.
I leave it to the reader to check that this category
is quasi-abelian.
A monomorphism in $\A$ is just an injective
morphism of sheaves. An epimorphism is a morphism of sheaves whose
cokernel is torsion. The kernel of a morphism of torsion-free
sheaves in $\A$ is just the usual sheaf-theoretic kernel, but the
cokernel in $\A$ is the usual cokernel modded out by its torsion
subsheaf. All epimorphisms are strict, whereas
a monomorphism is strict precisely
if its cokernel
as a map of sheaves is torsion-free.

\begin{lemma}
\label{beer}
Let $\P$ be a slicing of a triangulated category $\T$.
For any interval $I\subset \R$ of length $<1$,
the full subcategory $\P(I)\subset\T$ is
quasi-abelian. The strict short exact sequences in $\A$ are in
one-to-one correspondence with triangles in $\T$ all of whose vertices
are objects of $\T$.
\end{lemma}

\begin{pf}
Assume for definiteness that $I=(a,b)$ with $0<b-a<1$. The other cases are
equally easy. The result then follows by applying
 Lemma \ref{schn} to the embeddings
$\P((a,b))\subset\P((a,a+1])$ and $\P((a,b))\subset\P(([b-1,b))$ and
using Lemma \ref{blgl}.
\end{pf}

In what follows I shall abuse notation in a number of ways. Suppose
$A$, $B$ and $E$ are objects of a quasi-abelian category $\A$.
Then I shall
write $A\subset E$ to mean that there
is a strict monomorphism $i\colon A\to E$.
I shall also call $A$ a strict subobject of
$E$ and write $E/A$ for the cokernel of $i$.
Similarly, I write $E\onto B$ to mean that there is a strict
epimorphism $E\to B$ in $\A$ and refer to $B$ as a strict quotient of $E$.

As in the case of an abelian category, the
partial order $\subset$ allows one to say what it means for a
quasi-abelian category $\A$ to
be artinian or noetherian. Thus, for example, $\A$ is artinian if any infinite chain
\[\cdots\subset E_{j+1}\subset E_{j}\subset\cdots\subset E_2\subset
E_1\]
of strict subobjects in $\A$ must stabilise.
If $\A$ is artinian and noetherian then it
is said to be of finite length. For example, the category $\A$ of
torsion-free sheaves
described above is of finite length, because the rank function is additive
on $\A$ and every nonzero object of $\A$ has positive rank.

Using the notion of a strict subobject in a quasi-abelian
category, one can give a definition of semistability in a
quasi-abelian category, depending on a
choice of  stability function $Z\colon \K(\A)\to\C$.
Of course there is no reason to expect the resulting notion to have
good properties.
Nonetheless, the proof of Theorem \ref{lasty} 
will hinge on showing that in certain cases this notion of stability
in a quasi-abelian category
does in fact behave nearly as well as in the abelian case.

It will be convenient extend the definition so as to include possibly skewed stability functions as follows.

\begin{defn}
A \emph{skewed stability function} on
a quasi-abelian  category $\A$ is a group homomorphism $Z\colon \K(\A)\to\C$
such that there is a strict half-plane \[H_{\alpha}=\{r\exp(i\pi\phi):r>0\text{ and }\alpha<\phi\leq \alpha+1\}
\subset\C,\]
defined by some $\alpha\in\R$, such that
$Z(E)\in H_\alpha$ for all objects $0\neq E\in\A$.
\end{defn}

Clearly one can always reduce to the unskewed case $\alpha=0$ but
in fact it will
not always be convenient to do so.
Given a skewed stability function $Z\colon\K(\A)\to\C$, define
the phase of an object $0\neq E\in\A$ to be
\[\phi(E)=(1/\pi)\arg Z(E)\in (\alpha,\alpha+1].\]
An object $0\neq
E\in\A$ is then defined to be semistable if
for every strict subobject $0\neq A\subset E$ one has
$\phi(A)\leq\phi(E)$. An equivalent condition is that
$\phi(E)\leq\phi(B)$
for every nonzero strict quotient $E\onto B$.

A Harder-Narasimhan filtration of an object $0\neq E\in\A$ is a
finite chain of strict subobjects
\[0=E_0\subset E_1\subset \cdots\subset E_{n-1}\subset E_n=E\]
whose factors $F_j=E_j/E_{j-1}$ are semistable objects of $\A$ with
\[\phi(F_1)>\phi(F_2)>\cdots>\phi(F_n).\]
Recall that when $\A$ is abelian, Harder-Narasimhan
filtrations are
unique, essentially because if $f\colon E\to F$ is a nonzero map between
semistable objects
then $\phi(E)\leq\phi(F)$.
But the proof of
this fact depends on the assumption that all morphsims are strict,
so there is no reason to expect the corresponding result to hold in
the quasi-abelian context.

% ********************************************************************
% ********************************************************************
% ********************************************************************

\section{Stability conditions}
\label{stabilityconditions}

This section introduces the idea of a stability condition on a
triangulated category, which combines the notions of slicing and
stability function. The mathematical justification for this combination
seems to be that, as Theorem \ref{lasty} shows,
it leads to nice deformation properties.

\begin{defn}
A stability condition $\sigma=(Z,\P)$ on a triangulated category $\T$
consists of a group homomorphism
$Z\colon\K(\T)\to\C$ and a 
slicing $\P$ of $\T$
such that
if $0\neq E\in\P(\phi)$ then $Z(E)=m(E)\exp(i\pi \phi)$ for some $m(E)\in\R_{>0}$.
\end{defn}

The linear map $Z\colon\K(\T)\to\C$ will be referred to as the
\emph{central charge} of the stability condition.
The following Lemma shows that each category $\P(\phi)$ is
abelian.
The nonzero objects of $\P(\phi)$ are said to be
\emph{semistable}
in $\sigma$ of phase $\phi$, and the simple objects of $\P(\phi)$
are said to be \emph{stable}.

\begin{lemma}
If $\sigma=(Z,\P)$ is a stability condition on a triangulated
category $\T$ then each subcategory $\P(\phi)\subset\T$ is abelian.
\end{lemma}

\begin{pf}
The category $\P(\phi)$ is a full additive subcategory of the abelian category
$\A=\P((\phi-1,\phi])$. It will therefore be enough to show that if
$f\colon E\to F$ is
a morphism in $\P(\phi)$ then the kernel and cokernel of $f$,
considered as a morphism of $\A$, actually lie in $\P(\phi)$. But if
\[0\lra A\lra E\lra B\lra 0\]
is a short exact sequence in $\A$ and $E$ is an object of $\P(\phi)$
then
Lemma
\ref{blgl} implies that $B\in\P(\phi)$ and drawing a picture one sees
that $A\in\P(\phi)$ also.
\end{pf}

Let $\sigma=(Z,\P)$ be a stability condition on a triangulated
category $\T$.
Recall that
the decomposition of an object $0\neq E\in\T$ given in the definition
of a slicing is unique;
the objects $A_j$ will be called the \emph{semistable factors} of $E$
with respect to $\sigma$.
I shall write $\phi^{\pm}_{\sigma}(E)$ for $\phi^{\pm}_{\P}(E)$;
thus $\phi^+_{\sigma}(E)\geq\phi^-_{\sigma}(E)$ with equality
precisely if $E$ is semistable in $\sigma$.
The \emph{mass} of $E$ is defined to be
the positive real number
$m_{\sigma}(E)=\sum_i |Z(A_i)|$.
By the triangle inequality
one has $m_{\sigma}(E)\geq |Z(E)|$.
When the stability condition $\sigma$ is clear from the context
I often
drop it from the
notation and write $\phi^{\pm}(E)$ and $m(E)$.

The following result shows the relationship between t-structures and stability conditions.

\begin{prop}
\label{pg}
To give a stability condition on a triangulated category $\T$ is
equivalent to giving a bounded t-structure on $\T$ and a 
stability function on its heart with
the Harder-Narasimhan property.
\end{prop}

\begin{pf}
Note first that if $\A$ is the heart of a bounded t-structure on $\T$
then $\K(\A)$ can be identified with $\K(\T)$.
If $\sigma=(Z,\P)$ is a stability condition on
$\T$, the t-structure $\P(> 0)$ is bounded with
heart $\A=\P((0,1])$. The central charge $Z$ defines
a  stability function on $\A$ and it is easy to check
that the corresponding semistable objects are
precisely the nonzero objects of the categories $\P(\phi)$ for $0<\phi\leq 1$.
The decompositions of objects of $\A$ given by Definition
\ref{slicing}(c) are Harder-Narasimhan filtrations.

For the converse, suppose $\A$ is the heart of a bounded t-structure
on $\T$ and $Z\colon\K(\A)\to\C$ is a  stability function on $\A$
with the Harder-Narasimhan property. Define a stability condition
$\sigma=(Z,\P)$ on $\T$ as follows.
For each $\phi\in(0,1]$ let $\P(\phi)$ be the
full additive subcategory of $\T$
consisting of
semistable objects
of $\A$ with phase $\phi$, together with the zero objects of $\T$. 
The first condition of Definition \ref{slicing} then
determines $\P(\phi)$ for all $\phi\in\R$ and condition (b) is easily
verified. For any nonzero $E\in\T$ a filtration as in
Definition \ref{slicing}(c)
can be obtained by combining the decompositions of Lemma \ref{ll} with the
Harder-Narasimhan filtrations of nonzero objects of $\A$.
\end{pf}

I shall now give some examples of stability conditions.

\begin{example}
\label{standard}
Let $\A$ be the
category of coherent $\OO_X$-modules on a nonsingular projective
curve $X$ over an algebraically closed field $k$ of characteristic
zero, and set
$Z(E)=-\deg(E)+i\rk(E)$ as in the introduction.
Applying Proposition \ref{pg} gives a stability condition on the bounded
derived category $\D(\A)$.
\end{example}

This example will be considered in more detail in  Section
\ref{curve} below, where I study the set of all stability conditions
on the derived category of an elliptic curve.

\begin{example}
\label{quiver}
Let $A$ a finite-dimensional algebra over a field $k$.
Let $\A$ be the abelian category of finite-dimensional
left $A$-modules. Thus $\A$ is
a finite-length category whose Grothendieck group
$\K(\A)$ is isomorphic to the free abelian group on the
(finite) set of simple $A$-modules.
There is a group homomorphism
$r\colon\K(\A)\to\Z$
defined
by sending an $A$-module to its dimension as a vector space over $k$.
For any homomorphism $\lambda\colon \K(\A)\to\QQ$
the formula $Z(E)=\lambda(E)+i r(E)$ defines a slope function on
$\A$, and
Proposition \ref{pg} shows that each of these
slope functions determines a stability condition on the bounded
derived category $\D(\A)$.
\end{example}

The final example is rather degenerate and
is included purely to motivate
the introduction of the local-finiteness condition below.

\begin{example}
\label{bad}
Let $\A$ be the category of coherent $\OO_X$-modules on a nonsingular
projective
curve $X$ as in Example \ref{standard}, and let $(Z,\P)$ be the stability
condition on $\D(\A)$ defined there.
Let $0<\alpha<1/2$ be such that $\zeta=\tan
(\pi\alpha)$ is irrational. Then the
bounded t-structure $\P(>\alpha)=\P(\geq \alpha)\subset\T$ has heart
$\B=\P((\alpha,\alpha+1))$. Define a  stability function on $\B$ by the
formula
\[W(E)=i(\rk(E)+\zeta\deg(E)).\]
Note that all nonzero objects of $\B$ are semistable with the same
phase. Applying Proposition \ref{pg} gives a stability condition
$(W,\Q)$ on $\T$ such that $\Q(\ha)=\B$, and $\Q(\psi)=0$ unless
$\psi-\ha\in\Z$.
\end{example}

In order to eliminate such examples and to prove nice theorems
it will be useful to impose the following
extra condition on stability conditions.

\begin{defn}
\label{lofo}
A slicing $\P$ of a triangulated category $\T$ is \emph{locally-finite} if
there exists a real number $\eta>0$ such that for all $t\in\R$ the
quasi-abelian category $\P((t-\eta,t+\eta))\subset\T$ is of finite
length.
A stability condition $(Z,\P)$ is locally-finite if the corresponding
slicing $\P$ is.
\end{defn}

Note that if $(Z,\P)$ is a stability condition with a central charge
$Z$ whose image is a discrete subgroup of $\C$,
then for any interval $I\subset\R$
of length $<1$ the quasi-abelian category $\P(I)$ is of finite
length. Thus
 the first two examples of stability conditions
given above are
locally-finite. But the stability condition described in Example
\ref{bad} is not locally-finite in general, because as one can easily
check, the
abelian category $\B$ is not always of finite length.

% *************************************************************************
% *************************************************************************

\section{The space of stability conditions}
\label{space}

Fix a triangulated category $\T$
and write $\Slice(\T)$
for the set of locally-finite slicings
of $\T$ and $\Stab(\T)$ for the set of locally-finite stability conditions
on $\T$. The aim of this section is to define natural topologies on these
spaces. In fact, everything in this section applies equally well
without the locally-finite condition, which will only become important
in Section 7.

The first observation to be made is that the function
\[d(\P,\Q)=\sup_{0\neq E\in\T}
\bigg\{|\phi^-_{\P}(E)-\phi^-_{\Q}(E)|,|\phi^+_{\P}(E)-\phi^+_{\Q}(E)|
\bigg\}
\in[0,\infty]\]
defines a generalised metric\footnote{See the Notation section.} 
on $\Slice(\T)$. To check this one just needs to
note that if $d(\P,\Q)=0$ then
every nonzero object of $\P(\phi)$ is also an object of $\Q(\phi)$ so
that $\P=\Q$.
The following Lemma gives another way of writing this metric.

\begin{lemma}
If $\P$ and $\Q$ are slicings of a triangulated category $\T$ then
\[d(\P,\Q)=\inf\big\{\epsilon\in\R_{\geq 0}:
\Q(\phi)\subset\P([\phi-\epsilon,\phi+\epsilon])\text{ for all }\phi
\in\R\big\}\]
\end{lemma}

\begin{pf}
Write $d'(\P,\Q)$ for the expression in the statement of the
Lemma. First note that if $d(\P,\Q)\leq \epsilon$ then
for any nonzero $E\in\Q(\phi)$ one has
$\phi^+(E)\leq \phi+\epsilon$ and similarly $\phi^-(E)\geq
\phi-\epsilon$.
This implies that $\Q(\phi)\subset\P([\phi-\epsilon,\phi+\epsilon])$ and so
$d'(\P,\Q)\leq\epsilon$.

For the reverse inequality 
suppose $d'(\P,\Q)\leq \epsilon$ and take a nonzero object
$E\in\T$. Clearly if $E\in\Q(\leq \psi)$ then
$E\in\P(\leq \psi+\epsilon)$. But in the other direction, if
$E\notin\Q(\leq \psi)$ then there is some object
$A\in \Q(\phi)$ with $\phi>\psi$ and a nonzero map $A\to E$.
Since $\Q(\phi)\subset\P([\phi-\epsilon,\phi+\epsilon])$
it follows that $E\notin\P(\leq \psi-\epsilon)$.

These arguments show that
$|\phi^+_{\P}(E)-\phi^+_{\Q}(E)|\leq \epsilon$, and a similar argument
with $\phi^-$ completes the proof that $d(\P,\Q)\leq\epsilon$.
\end{pf}

Consider the inclusion of sets
\[\Stab(\T)\subset\Slice(\T)\times\Hom_{\Z}(\K(\T),\C).\]
When $\K(\T)$ has finite rank, one can give the vector space on the
right the standard topology, and obtain an
induced topology on $\Stab(\T)$. In general however, one has
to be a little careful, since there is no obviously natural choice of
topology on $\Hom_{\Z}(\K(\T),\C)$.

For each $\sigma=(Z,\P)\in\Stab(\T)$, define a function
\[ \|\cdot\|_{\sigma}\colon \Hom_{\Z}(\K(\T),\C)\to [0,\infty]\]
by sending a linear map $U\colon \K(\T)\tensor\C\to\C$ to
\[\|U\|_{\sigma}=\sup\bigg\{\frac{|U(E)|}{|Z(E)|}:\text{$E$
semistable in $\sigma$}\bigg\}.\]
Note that $\|\cdot\|_{\sigma}$ has all the properties of a norm
on the complex vector space $\Hom_{\Z}(\K(\T),\C)$, except that it may not be
finite.

For each real number $\epsilon\in(0,\quarter)$, define a subset
\[B_{\epsilon}(\sigma)=\{\tau=(W,\Q):\|W-Z\|_{\sigma}<\sin(\pi\epsilon)\text{
and }d(\P,\Q)<\epsilon\}\subset\Stab(\T).\]
To understand this definition note that the
condition $\|W-Z\|_{\sigma}<\sin(\pi\epsilon)$ implies
that for all objects $E$ semistable in $\sigma$, the phase of $W(E)$
differs from the phase of $Z(E)$ by less than $\epsilon$.

I claim that as $\sigma$ varies in $\Stab(\T)$ the subsets
$B_{\epsilon}(\sigma)$ form a
basis for a topology on $\Stab(\T)$. This boils down to the statement
that if $\tau\in B_{\epsilon}(\sigma)$ then
there is an $\eta>0$ such that $B_{\eta}(\tau)\subset
B_{\epsilon}(\sigma)$, which follows easily from the following crucial lemma.

\begin{lemma}
\label{pent}
If $\tau=(W,\Q)\in B_{\epsilon}(\sigma)$ then there are constants $k_i>0$
such that
\[ k_1\|U\|_{\sigma}< \|U\|_{\tau} <k_2\|U\|_{\sigma}\]
for all $U\in\Hom_{\Z}(\K(\T),\C)$.
\end{lemma}

\begin{pf}
First, note that for any stability condition $\sigma=(Z,\P)$ on $\T$,
and any real number $0\leq\eta<\half$, one has
\begin{equation*}
\tag{$*$}|U(E)|<\frac{\|U\|_{\sigma}}{\cos(\pi\eta)}|Z(E)|,
\end{equation*}
for every $0\neq E\in\T$ satisfying
$\phi^+_{\sigma}(E)-\phi^-_{\sigma}(E)<\eta$,
and for all linear maps $U\colon \K(\T)\to\C$. To see this, just
decompose
$E$ into semistable factors $A_1,\cdots ,A_n$ in $\sigma$, apply
the definition of $\|U\|_{\sigma}$ to each object $A_i$, and note that
the points $Z(A_i)\in\C$ lie in a sector bounded by an angle of at
most $\pi\eta$.

Now consider the situation of the Lemma.
Since $d(\P,\Q)<\quarter$ and $\|Z-W\|_{\sigma}<\infty$,
one can apply $(*)$ with $U=W-Z$
to show that there is a constant $\kappa>0$ with
$|Z(E)|<\kappa |W(E)|$ for all $E$ semistable in $\tau$.
Take a linear map $U\colon\K(\T)\to\C$ and
suppose $E$ is semistable in $\tau$. Since $d(\P,\Q)<\epsilon$
it follows that $(*)$ holds with $\eta=2\epsilon$. Combining this with
the above inequality gives $\|U\|_{\tau}<k_2\|U\|_{\sigma}$.
But now $\|Z-W\|_{\tau}<\infty$ so one can
swap $\sigma$ and $\tau$ and
repeat the argument to give the reverse inequality.
\end{pf}

Equip $\Stab(\T)$ with the topology generated by the
basis of open sets $B_{\epsilon}(\sigma)$.
Let $\Sigma$ be a connected component of
$\Stab(\T)$. By
Lemma \ref{pent}, the subspace
\[\{U\in\Hom_{\Z}(\K(\T),\C):\|U\|_{\sigma}<\infty\}\subset
\Hom_{\Z}(\K(\T),\C)\]
is locally constant on $\Stab(\T)$ and hence constant on
$\Sigma$. Denote it by $V(\Sigma)$. Note
that if $\sigma=(Z,\P)\in\Sigma$ then $Z\in V(\Sigma)$. 
Note also
that for each $\sigma\in\Sigma$ the function $\|\cdot\|_{\sigma}$
defines a norm on $V(\Sigma)$, and that by Lemma \ref{pent}, all these
norms are equivalent. Thus one has
 
\begin{prop}
\label{p}
For each connected component $\Sigma\subset\Stab(\T)$
there is a linear subspace $V(\Sigma)\subset \Hom_{\Z}(\K(\T),\C)$
with a well-defined linear topology and a continuous map
$\ZZ\colon\Sigma\to V(\Sigma)$ which sends a stability condition $(Z,\P)$
to its central charge $Z$.\qed
\end{prop}

The proof of Theorem \ref{lasty} will be completed in Section \ref{deform}
by showing
that the map
$\ZZ$ of Proposition \ref{p} is a local homeomorphism.
The
following lemma shows that $\ZZ$ is at least locally injective.

\begin{lemma}
\label{inj}
Suppose $\sigma=(Z,\P)$ and $\tau=(Z,\Q)$ are
stability conditions on $\T$ with
the same central
charge $Z$. Suppose also that $d(\P,\Q)< 1$.
Then $\sigma=\tau$.
\end{lemma}

\begin{pf}
Suppose on the contrary that $\sigma\neq \tau$.
Then there is a nonzero
object $E\in\P(\phi)$ which is not an element of $\Q(\phi)$.
One could not have $E\in\Q(\geq\phi)$ because the assumption that
$d(\sigma_1,\sigma_2)<1$ would then imply that $E\in\Q([\phi,\phi+1))$
which contradicts the fact that $\sigma$ and $\tau$ have the same
central charge. Similarly one could not have $E\in\Q(\leq\phi)$.
Thus there is a triangle
\[\xymatrix@C=.5em{A \ar[rrrr] &&&& E  \ar[dll] \\
 && B \ar@{-->}[ull]} \]
with $A\in\Q((\phi,\phi+1))$ and $B\in\Q((\phi-1,\phi])$ nonzero.
One cannot have $A\in\P(\leq \phi)$
because this would imply $A\in\P((\phi-1,\phi])$ contradicting the
fact that $\sigma$ and $\tau$ have the same central charge.
Thus
there is an object $C\in\P(\psi)$ with $\psi>\phi$ and a nonzero morphism
$f\colon C\to A$. The composite map $C\to E$ must be zero so $f$
factors via $B[-1]$. Since $B[-1]\in\Q(\leq\phi-1)$
this gives a contradiction.
\end{pf}

% *************************************************************************
% *************************************************************************
% *************************************************************************
% *************************************************************************

\section{Deformations of stability conditions}
\label{deform}

In this section I complete the proof of Theorem \ref{lasty} by proving
a result that allows one to lift deformations of the central charge
$Z$ to deformations of stability conditions. It was Douglas' work
that first
suggested that such a result might be true.

\begin{thm}
\label{biggy}
Let $\sigma=(Z,\P)$ be a locally-finite
stability condition on a triangulated category $\T$. Then there is an
$\epsilon_0>0$ such that if $0<\epsilon<\epsilon_0$ and
$W\colon\K(\T)\to\C$ is a group homomorphism satisfying
\[|W(E)-Z(E)|<\sin(\pi\epsilon)|Z(E)|\]
for all $E\in\T$ semistable in $\sigma$,
then there is a locally-finite stability condition
$\tau=(W,\Q)$ on $\T$ with $d(\P,\Q)<\epsilon$.
\end{thm}

After what was proved in Section \ref{space}
this will be enough to
give Theorem \ref{lasty}.
Note that Lemma \ref{inj} shows that, providing $\epsilon_0<1/2$,
the stability condition $\tau$ of Theorem
\ref{biggy} is unique. The reader should
think of the number $\epsilon_0$ as being very small. In fact,
it will be enough to
assume that $\epsilon_0<1/8$ and
that each of the quasi-abelian categories
$\P((t-4\epsilon_0,t+4\epsilon_0))$ has finite length. Since
$\Q((t-\epsilon,t+\epsilon))\subset \P((t-2\epsilon,t+2\epsilon))$
for all $t$, the condition that $\tau$ be locally-finite is automatic.
The proof of the theorem will be broken up into a series of lemmas.
Throughout notation
will be fixed as in the statement of the Theorem. In particular,
$W\colon\K(\T)\to\C$ is a group homomorphism satisfying the hypotheses
of the Theorem, and $0<\epsilon<\epsilon_0$ is a fixed real number.

\begin{defn}
A \emph{thin subcategory} of $\T$ is a full subcategory
of the form $\P((a,b))\subset \T$ where $a$ and $b$ are real numbers
with $0<b-a<1-2\epsilon$.
\end{defn}

Note that any thin subcategory of $\T$ is quasi-abelian.
Recall that the condition on $W$ in the statement of the theorem
implies that if $E$ is semistable in $\sigma$, then the phases of the
points $W(E)$ and $Z(E)$ differ by at most $\epsilon$.
It follows that 
if $\A=\P((a,b))$ is thin then $W$
defines a skewed stability function on $\A$.
To avoid confusion, the objects of $\A$ which are
semistable with respect to this stability function will be called $W$-semistable.
Also, given a nonzero object $E\in\A$,
write $\phi(E)$ for the phase of $Z(E)$ lying in the interval
$(a,b)$, and $\psi(E)$ for
the phase of $W(E)$ lying in the interval $(a-\epsilon, b+\epsilon)$. 

\begin{lemma}
\label{peasy}
Suppose $E$ is $W$-semistable
in some thin subcategory $\A\subset \T$, and set $\psi=\psi(E)$.
Then $E\in\P((\psi-\epsilon,\psi+\epsilon))$.
\end{lemma}

\begin{pf}
Put $\phi=\phi^+(E)$. There is a strict short exact sequence
\[0\lra A\lra E\lra B\lra 0\]
in $\A$ such that $A\in\P(\phi)$ and $B\in\P(<\phi)$. Then 
$\psi(A)\leq\psi(E)$ because $E$ is $W$-semistable. But as above, one has
$\psi(A)\in(\phi-\epsilon,\phi+\epsilon)$ and so it follows that
$\phi<\psi+\epsilon$.
A similar argument shows that  $\phi^-(E)>\psi-\epsilon$.
\end{pf}

This notion of $W$-semistability for an object $E$ of a
thin subcategory
is too weak unless $E$ lies
well inside $\A$ in a certain sense. The problem is that if
$E$ lies near the boundary of $\A$ then there are not enough objects in $\A$
to destabilise $E$. This prompts the following definition.

\begin{defn}
Suppose $\A=\P((a,b))$ is a thin subcategory of $\T$.
A nonzero object $E\in\A$ is said to be \emph{enveloped} by
$\A$ if $a+\epsilon\leq\psi(E)\leq b-\epsilon$.
\end{defn}

The next Lemma shows that
with this idea one gets a notion of semistability which is independent
of a particular choice of thin subcategory.

\begin{lemma}
Suppose an object $E\in\T$ is enveloped by thin subcategories
$\B$ and $\CC$ of $\T$. Then $E$ is $W$-semistable in $\B$
precisely if it is $W$-semistable in $\CC$.
\end{lemma}

\begin{pf}
After Lemma \ref{peasy} one may as well assume that
$E$ is enveloped by the thin subcategory
$\P((\psi(E)-\epsilon,\psi(E)+\epsilon))$.
Thus it is enough to treat the case when
$\B\subset\CC$, and in fact, by the symmetry
of the situation, one can assume that $\B=\P((a,b))$ and
$\CC=\P((a,c))$ for real numbers $a<b<c$. Of course,
if $E$ is $W$-semistable in $\CC$ then it
is also $W$-semistable in $\B$, because any strict
short exact sequence in $\B$
is also a strict
short exact sequence in $\CC$.

For the converse, suppose $E$
is unstable in $\CC$ so that there is a strict short exact sequence in $\CC$
\[0\lra A\lra E\lra B\lra 0\]
with $\psi(A)>\psi(E)>\psi(B)$. Then, by Lemma \ref{blgl}, one has
$\phi^+(A)\leq\phi^+(E)$, so  since $E\in\B$, one has $A\in\B$ also.
There is a strict short exact sequence
\[0\lra B_1\lra B\lra B_2\lra 0\]
in $\CC$ with $B_1\in\P([b,c))$ and $B_2\in\B$. Note that because $E$
is enveloped by $\B$ one has $\psi(E)\leq b-\epsilon<\psi(B_1)$.
Consider the commuting diagram of strict short exact
sequences in $\CC$
\[
\begin{CD}
@. 0 @. 0 @. @.\\
@. @VVV @VVV @.@.\\
@. A @= A @.@. \\
@. @VVV @VVV @.@.\\
0 @>>> K @>>> E @>>> B_2 @>>> 0 \\
@. @VVV @VVV @| @. \\
0 @>>> B_1 @>>> B @>>> B_2 @>>> 0\\
@. @VVV @VVV @.@.\\
@.  0 @. 0 @. @.
\end{CD}\]

Then by Lemma \ref{blgl} again, $\phi^+(K)\leq\phi^+(E)$ and hence
$0\to K\to E\to B_2\to 0$
is a strict short exact
sequence in $\B$. But $\psi(K)>\psi(E)$ and therefore
$E$ is not
$W$-semistable in $\B$.
\end{pf}

For each
$\psi\in\R$ define $\Q(\psi)\subset\T$ to be the
full additive subcategory of $\T$ consisting of the zero objects of
$\T$ together with those objects
$E\in\T$ which are
$W$-semistable of phase $\psi$ in some thin enveloping subcategory
$\P((a,b))$. To prove Theorem
\ref{biggy} it must be shown that the pair $(W,\Q)$ defines a stability
condition on $\T$. The following Lemma gives axiom (c) of
Definition \ref{pemb}.

\begin{lemma}
\label{fa}
If $E\in\Q(\psi_1)$ and $F\in\Q(\psi_2)$ and $\psi_1>\psi_2$ then
$\Hom_{\T}(E,F)=0$.
\end{lemma}

\begin{pf}
Suppose instead that there is a nonzero map $f\colon E\to F$.
By Lemma \ref{peasy} this implies that $\psi_1-\psi_2<2\epsilon$.
Set $a=(\psi_1+\psi_2)/2-1/2$
and consider the abelian subcategory
$\A=\P((a,a+1])\subset\T$ which contains $E$ and $F$. In the abelian
category $\A$ there are short
exact sequences
\[0\lra\ker f\lra E\lra \im f\lra 0\]
and
\[0\lra \im f \lra F\lra \coker f \lra 0.\]
By Lemmas \ref{blgl} and \ref{peasy},
one has $\ker f\in\P((a,\psi_1+\epsilon))$, $\coker
f\in\P((\psi_2-\epsilon, a+1])$ and $\im
f\in\P((\psi_1-\epsilon,\psi_2+\epsilon))$. Providing $\epsilon$
is small enough (say $\epsilon<1/8$),
there is a thin subcategory of $\T$ enveloping $E$
in which the first
sequence is strict short exact, and similarly a thin subcategory
enveloping  $F$
in which the second sequence is strict short exact. Since $E$ and $F$ are
$W$-semistable in any enveloping category it follows that
$\psi_1\leq\psi(\im f)\leq\psi_2$, a contradiction.
\end{pf}

The last step in the proof of Theorem \ref{biggy} is to construct
filtrations of objects of $\T$ with factors in
the subcategories $\Q(\psi)$.

\begin{lemma}
\label{floor}
Let $\A=\P((a,b))\subset \T$ be a thin subcategory of finite length.
Then every nonzero object of $\P((a+2\epsilon,b-4\epsilon))$ has a
finite Harder-Narasimhan filtration whose factors 
are $W$-semistable objects of $\A$ which are enveloped by $\A$.
\end{lemma}

\begin{pf}
The proof goes along the same lines as that of
Proposition \ref{rud}, replacing subobjects by strict subobjects and
quotients by strict quotients. Here I just indicate the necessary
changes.
Clearly the chain conditions hold because of the assumption that $\A$
has finite length.
Note also that if an object $E\in\P((a+2\epsilon,b-4\epsilon))$ has a
Harder-Narasimhan filtration with $W$-semistable
factors $F_1,\cdots, F_n$, then
$\psi(F_1)\geq\psi(E)>a+\epsilon$, and the
fact that there is a nonzero map $F_1\to E$ together with Lemma
\ref{peasy} ensures that $\psi(F_1)<b-3\epsilon$. In this way one sees
that the factors of $E$ are automatically enveloped by $\A$.

Define $\mathcal G$ be the class of of nonzero objects
$E\in\P((a,b-4\epsilon))$
for which every nonzero strict quotient $E\onto B$ in $\A$ satisfies
$\psi(B)>a+\epsilon$. By Lemma \ref{blgl} the class $\mathcal G$
contains all nonzero objects of $\P((a+2\epsilon,b-4\epsilon))$,
so it will be enough to show that all objects of
$\mathcal G$ have a Harder-Narasimhan filtration. The commutative
diagram $(\dagger)$ and Lemma \ref{blgl} show that if
\[0\lra E'\lra E\lra B\lra 0\]
is a strict short exact sequence in $\A$ with $E\onto B$ a maximally
destabilising quotient (mdq) and
$E\in\mathcal G$ then $E'\in\mathcal G$ also. Thus the inductive step
in the proof of Proposition \ref{rud} stays within the class $\G$ and
it will be enough
to show that every object in $\mathcal G$ has a mdq.

To make the induction work it is helpful to prove the existence of maximally
destabilising quotients for a larger
class of objects $\mathcal H$, namely nonzero objects $E\in\A$ with
$\psi(E)<b-3\epsilon$ such that every nonzero strict quotient $E\onto B$ in
$\A$ satisfies $\psi(B)>a+\epsilon$. Note that if $E\in\mathcal H$ and
$E\onto E'$ is a nonzero strict quotient with $\psi(E)\geq\psi(E')$ then
$E'\in\mathcal H$ also.

Suppose then that $E\in\mathcal H$.
The key observation is that I can always assume that
$\phi^+(E)<\psi(E)+\epsilon$. Otherwise there is a strict short exact
sequence
\[0\lra A\lra E\lra E'\lra 0\]
with $A\in\P(\geq \psi(E)+\epsilon)$ and
$E'\in\P(<\psi(E)+\epsilon)$. Note that $\psi(A)>\psi(E)>\psi(E')$.
I claim that if $E'\onto B$ is a maximally
destabilising quotient for
$E'$ then the composite map $E\onto B$ is a maximally
destabilising quotient for $E$. Indeed, if
$E\onto B'$ is a $W$-semistable quotient in $\A$ with $\psi(B')\leq\psi(B)$
then $\psi(B')\leq\psi(E)$ and so by Lemma \ref{peasy}
one has $\phi^+(B')<\psi(E)+\epsilon$. It follows that
$\Hom_{\A}(A,B')=0$ and hence $E\onto B'$ factors via $E'$. This
proves the claim.

By Lemmas \ref{blgl} and
\ref{peasy}, the inequalities $\phi^+(E)<\psi(E)+\epsilon$ and
$\psi(E)<b-3\epsilon$ are
enough to guarantee 
that every
$W$-semistable strict subobject of $E$ is enveloped by $\A$. By definition
of the class $\mathcal H$,
every $W$-semistable strict quotient of $E$ is
also enveloped by $\A$. Thus, using Lemma \ref{fa},
the argument of Proposition \ref{rud} can be
applied as in the abelian case to show that $E$ has a maximally
destabilising quotient.
\end{pf}

For each real number $t$ define $\Q(>t)$ to
be the full extension-closed  subcategory of
$\T$ generated by the subcategories $\Q(\psi)$ for $\psi>t$. Similarly
define full subcategories $\Q(\leq t)\subset\T$ and $\Q(<t)\subset
\T$.

I claim that $\Q(>t)$
is a t-structure on $\T$.
To prove this I must show that for every
$E\in\T$ there is a triangle
\[A\lra E\lra B\]
 with $A\in\Q(>t)$ and
$B\in\Q(\leq t)$. But note that Lemmas \ref{peasy} and \ref{floor} show that
$\P(s)$ is contained in the subcategory $\Q(>t)$ for
$s\geq t+\epsilon$ and in the subcategory $\Q(< t)$ for $s\leq t-\epsilon$.
Thus it will be enough to
consider the case when $E\in\P((t-\epsilon,t+\epsilon))$. Consider $E$ as an
object of the quasi-abelian category $\P((t-3\epsilon,t+5\epsilon))$
which has finite length by the assumptions on $\epsilon$. Applying
Lemma \ref{floor} gives a Harder-Narasimhan filtration of $E$ which is
enough to prove the claim.

The final step in the proof of Theorem \ref{biggy}
is to show that every nonzero object of $\T$
has a finite filtration by
objects of the subcategories $\Q(\psi)$. It will be enough to prove
this for objects in each of the full subcategories
\[\Q((t,t+\delta))=\Q(>t)\cap\Q(< t+\delta)\] for some small $\delta>0$.
The result then follows by embedding $\Q((t,t+\delta))$ in the finite length
quasi-abelian subcategory
$\P((t-3\epsilon,t+5\epsilon+\delta))$ and applying Lemma \ref{floor}.

% *************************************************************************
% *************************************************************************

\section{More on the space of stability conditions}

This section contains a couple of general results about spaces
of stability conditions.  The first
shows that the topology on $\Stab(\T)$ defined in Section
\ref{space} can be induced by a
natural metric. Since this result is not necessary for
Theorem \ref{lasty}
some of the details of the proof are left to the reader.

\begin{prop}
\label{spiritair}
Let $\T$ be a triangulated category.
The function
\[d(\sigma_1,\sigma_2)=\sup_{0\neq E\in\T}
\bigg\{|\phi^-_{\sigma_2}(E)-\phi^-_{\sigma_1}(E)|,|\phi^+_{\sigma_2}(E)-\phi^+_{\sigma_1}(E)|
,|\log \frac{m_{\sigma_2}(E)}{m_{\sigma_1}(E)}|\bigg\}
\in[0,\infty]\]
defines a generalised metric on $\Stab(\T)$.
The induced topology is the same as that defined in Section
\ref{space}.
\end{prop}

\begin{pf}
It is easy to see that the given formula defines a generalised metric,
the only thing to check is that if $d(\sigma_1,\sigma_2)=0$ then
$\sigma_1=\sigma_2$. But $d(\sigma_1,\sigma_2)=0$ implies that an
object $E\in\T$ is semistable in $\sigma_1$ precisely if it is
semistable in $\sigma_2$, and that for any nonzero $E$ one has
$m_{\sigma_1}(E)=m_{\sigma_2}(E)$. It follows that the central charges
of $\sigma_1$ and $\sigma_2$ are the same, since they agree on
semistables and these span the Grothendieck group $\K(\T)$.

To prove that the topology induced by $d(-,-)$ is the same as the one
given by the basis of open sets $B_{\epsilon}(\sigma)$
one must first show that the sets $B_{\epsilon}(\sigma)$ are
open in the topology induced by the metric.
This boils down to the statement that if
$\tau=(W,\Q)\in\Stab(\T)$ is a small enough distance from
$\sigma=(Z,\P)$ then \[|W(E)-Z(E)|<\sin(\pi\epsilon)|Z(E)|\] for all
objects $E\in\T$ semistable in $\sigma$. This is easy enough to see and is
probably best done privately
with a picture.

The reverse implication requires a little more care. Take
$\sigma=(Z,\P)\in\Stab(\T)$ and fix a constant $\kappa>1$.
What one needs to show
is that for small enough $\epsilon>0$ the set $B_{\epsilon}(\sigma)$
has the property that
\[\tau=(W,\Q)\in B_{\epsilon}(\sigma)
\implies m_{\tau}(E)<\kappa m_{\sigma}(E)\text{ for all }0\neq E\in\T.\]

Suppose then that $\tau\in B_{\epsilon}(\sigma)$ and
consider first the case when
$\phi^+_{\sigma}(E)-\phi^-_{\sigma}(E)<\eta$
for some $\eta\in(0,\half)$. Split $E$ into semistable factors $A_i$
with respect to $\tau$. Then
$\phi^+_{\sigma}(A_i)-\phi^-_{\sigma}(A_i)<2\epsilon$ for each $i$, so that
equation $(*)$ of the proof of Lemma \ref{pent} gives
\[|W(A_i)|<\big(1+\frac{\sin(\pi\epsilon)}{\cos(2\pi\epsilon)}\big)|Z(A_i)|.\]
Since the vectors $Z(A_i)$ lie in a sector bounded by an angle of at most
$\pi(4\epsilon+\eta)$, and $|Z(E)|\leq m_{\sigma}(E)$,
it follows that there is a constant 
$r(\epsilon,\eta)>1$ such that
\[m_{\tau}(E)<r(\epsilon,\eta) m_{\sigma}(E),\]
and that moreover $r(\epsilon,\eta)\to 1$ as $\max(\epsilon,\eta)\to 0$.

Consider now a general nonzero object $E\in\T$.
Fix real numbers $\phi$ and a positive integer $n$. For each integer $k$,
define intervals
\[I_k=\big[\phi+kn\epsilon,\phi+(k+1)n\epsilon\big),\quad
J_k=\big[\phi+(kn-1)\epsilon,\phi+((k+1)n+1))\epsilon\big),\]
and let $\alpha_k$
and $\beta_k$ be the truncation functors projecting into the
subcategories
$\Q(I_k)$ and $\P(J_k)$ respectively. It is an easy enough exercise
to check that $d(\P,\Q)<\epsilon$ implies
$\alpha_k\circ \beta_k=\alpha_k$
so that
\[m_{\tau}(E)=\sum_k m_{\tau}(\alpha_k(E))\leq \sum_k
m_{\tau}(\beta_k(E))<r(\epsilon,(n+2)\epsilon)\sum_k
m_{\sigma}(\beta_k(E)).\]
But now one can choose $\phi$ so that
\[\sum_k m_{\sigma}(\beta_k(E))\leq
\bigg(1+\frac{2}{n}\bigg)m_{\sigma}(E),\]
so that if one lets $\epsilon\to 0$ and $n\to\infty$ in such a way
that $n\epsilon\to 0$ then one sees that for small enough
$\epsilon$ one has $m_{\tau}(E)<\kappa
m_{\sigma}(E)$ for all nonzero $E\in\T$.
\end{pf}

Proposition \ref{spiritair}
has the consequence that for any nonzero object $E\in\T$
the functions
\[\phi^{\pm}(E):\Stab(\T)\to \R\quad\mbox{ and }\quad
m(E):\Stab(\T)\to\R_{>0}\]
are continuous.
It follows immediately from this that the subset of
$\Stab(\T)$ consisting of those stability conditions in which a given object
$E\in\T$ is semistable is a closed subset. Indeed, if $E$ is nonzero,
it is precisely the set
of $\sigma\in\Stab(\T)$ for which the equality
$\phi^+_{\sigma}(E)=\phi^-_{\sigma}(E)$ holds.

\begin{lemma}
\label{groupactions}
The generalised metric space $\Stab(\T)$ carries a right action of the
group
$\grp$, the universal
covering space of $\GL(2,\R)$, and a left action by isometries of the 
group
$\Aut(\T)$ of exact autoequivalences of $\T$. These two actions
commute.
\end{lemma}

\begin{pf}
First note that the group $\grp$ can be thought of as the
set of pairs $(T,f)$ where $f\colon \R\to\R$ is an
increasing map with $f(\phi+1)=f(\phi)+1$, and  $T\colon \R^2\to\R^2$
is an orientation-preserving linear isomorphism, such that the induced maps on
$S^1=\R/2\Z=\R^2/\R_{>0}$ are the same.

Given a stability condition $\sigma=(Z,\P)\in\Stab(\T)$, and a
pair $(T,f)\in\grp$, define a new stability
condition $\sigma'=(Z',\P')$ by setting $Z'=T^{-1}\circ Z$ and
$\P'(\phi)=\P(f(\phi))$. Note that the semistable objects of the
stability conditions $\sigma$
and $\sigma'$ are the same, but the phases have been
relabelled.

For the second action, note that an element
$\Phi\in\Aut(\T)$ induces an automorphism $\phi$ of $\K(\T)$. If
$\sigma=(Z,\P)$ is a stability condition on $\T$ define $\Phi(\sigma)$
to be the stability condition $(Z\circ\phi^{-1},\P')$, where
$\P'(t)$=$\Phi(\P(t))$. The reader can easily check that this action
is by isometries and commutes with the first.
\end{pf}

It might be said that the existence of a $\grp$ action on $\Stab(\T)$
means that stability functions $Z$ should be considered as maps to $\R^2$
rather than maps to $\C$. At present I have no convincing
argument against this.

% *************************************************************************
% *************************************************************************

\section{Stability conditions on curves}

\label{curve}
Let $X$ be a nonsingular projective curve of genus one over $\C$, and let $\D(X)$ denote the bounded
derived category of coherent $\OO_X$-modules.
As in the introduction, $\Stab(X)$ will denote the space of
locally-finite numerical stability conditions on $\D(X)$.

Set $\K(X)=\K(\D(X))$ and write $\N(X)$ for the numerical
Grothendieck group $\N(\D(X))$ defined in Section 1.3.
The Riemann-Roch theorem shows that $\N(X)$ can be identified with
$\Z\oplus\Z$, with the quotient map $\K(X)\to\N(X)$
sending a class $[E]\in\K(X)$ to the
pair consisting of its rank
and degree.
The Euler form on $\N(X)$ is then given by
\[\eu((r_1,d_1),(r_2,d_2))=r_1 d_2-r_2 d_1.\]
As in Example \ref{standard}, there is a stability condition
$\sigma=(Z,\P)\in\Stab(X)$ with \[Z(E)=-\deg(E)+i\rk(E),\]
in which the objects of
the subcategories $\P(\phi)$ consist of shifts of semistable sheaves
on $X$, and whose heart is the category of coherent
$\OO_X$-modules. It follows from Lemma \ref{groupactions} and Theorem
\ref{lasty}
that there is a local homeomorphism
\[\ZZ\colon\Stab(X)\to \Hom_{\Z}(\N(X),\C)\]
whose image is some open subset of the two-dimensional vector space $\Hom_{\Z}(\N(X),\C)$.

\begin{thm}
The action of the group
$\grp$ on $\Stab(X)$
is free and transitive, so that
\[\Stab(X)\isom\grp.\]
\end{thm}

\begin{pf}
First note that if $E$ is an indecomposable sheaf on
$X$ then $E$ must be semistable
in any stability condition $\sigma\in\Stab(X)$
because otherwise there is a
nontrivial triangle $A\to E\to B$ with $\Hom_{\D(X)}(A,B)=0$, and
then Serre duality gives
\[\Hom^1_{\D(X)}(B,A)=\Hom_{\D(X)}(A,B)^*=0,\]
which implies that $E$ is
a direct
sum $A\oplus B$.

Take an element $\sigma=(Z,\P)\in\Stab(X)$. Suppose for a
contradiction that
the image of the central charge
$Z$ is contained in a real line in
$\C$. Since $\sigma$ is locally-finite, the heart $\A$ of
$\sigma$ must then be of finite length. If $A$ and $B$ are simple objects
of $\A$ then \[\Hom_{\D(X)}(A,B)=\Hom_{\D(X)}(B,A)=0,\] and it follows
from this that
$\eu(A,B)=0$. But this implies that all simple objects of $\A$ lie on the
same line in $\N(X)$, and hence that all objects of $\D(X)$ do too, 
which gives a contradiction.
Thus $Z$, considered as a map from $\N(X)\tensor\R=\R^2$
to $\C\isom\R^2$ is an isomorphism, and it follows that the action of
$\grp$ on $\Stab(X)$ is free.

Suppose $A$ and $B$ are line bundles on $X$ with $\deg(A)<\deg(B)$.
Since $A$ and $B$ are indecomposable they are semistable in
$\sigma$ with phases $\phi$ and $\psi$ say. The existence of maps
$A\to B$ and $B\to A[1]$ gives inequalities $\phi\leq\psi\leq\phi+1$,
which implies that $Z$ is orientation preserving.
Thus acting by an element of $\grp$, one can assume that
$Z(E)=-\deg(E)+i\rk(E)$, and that for some point $x\in X$ the
skyscraper sheaf
$\OO_x$ has phase $1$. Then all
semistable vector bundles on $X$ are semistable in $\sigma$ with phase
in the interval $(0,1)$, and it follows quickly from this that
$\sigma$ is the standard stability condition
described in Example \ref{standard}.
\end{pf}

The quotient $\Stab(X)/\Aut
\D(X)$ is also of interest. One can easily show that the autoequivalences of $\D(X)$
are generated by shifts, automorphisms of $X$ and twists by line
bundles together with the Fourier-Mukai transform \cite{Mu1}. Automorphisms of $X$ and twists by
line bundles of degree zero
act trivially on $\Stab(X)$ and one obtains
\[{\Stab(X)}\left/{\Aut \D(X)}\right.\isom
{\GL(2,\R)}\left/ \SL(2,\Z).\right.\]
which is easily seen to be a $\C^*$-bundle over the moduli space
of elliptic curves.

% *************************************************************************
% *************************************************************************

\bigskip

\noindent Department of Pure Mathematics,
University of Sheffield,
Hicks Building, Hounsfield Road, Sheffield, S3 7RH, UK.

\smallskip

\noindent email: {\tt t.bridgeland@sheffield.ac.uk}


\begin{thebibliography}{16}

\bibitem{Do3} P. Aspinwall and M.R. Douglas,
D-Brane stability and monodromy, preprint hep-th/0110071.

\bibitem{Be} A. Beilinson,
Coherent sheaves on $\PP\sp{n}$ and problems in linear algebra, 
Funktsional. Anal. i Prilozhen. 12 (1978), no. 3, 68--69;
English transl. in Functional Anal. Appl. 12 (1978), no. 3.

\bibitem{BBD} A. Beilinson, J. Bernstein and P. Deligne, Faisceaux
Pervers, Ast{\'e}risque {100}, Soc. Math de France (1983).

\bibitem{Bo1} { A.I. Bondal and D.O. Orlov,} Reconstruction of a variety from the
derived category and groups of autoequivalences, Compositio Math. 125
(2001), no. 3, 327--344, also math.AG/9712029.

\bibitem{VdB} A.I. Bondal and M. Van den Bergh, Generators and
representability of functors in commutative and noncommutative
geometry, preprint math.AG/0204218.


\bibitem{Br} {T. Bridgeland,} Stability conditions on K3 surfaces,
in preparation.
 
\bibitem{Do1} M.R. Douglas, D-Branes on Calabi-Yau manifolds, preprint
math.AG/0009209.

\bibitem{Do2} M.R. Douglas, D-branes, categories and $N=1$
supersymmetry. Strings, branes, and M-theory. J. Math. Phys. 42
(2001), no. 7, 2818--2843, also hep-th/0011017.

\bibitem{Do4} M.R. Douglas, Dirichlet branes,
homological mirror symmetry, and stability, preprint math.AG/0207021,
to appear in the 2002 ICM proceedings.

\bibitem{GM} S.I. Gelfand and Yu. I. Manin, Methods of Homological
Algebra, Springer-Verlag (1996).


\bibitem{HN} G. Harder and M.S. Narasimhan,
On the cohomology groups of moduli spaces of vector bundles on curves,
Math. Ann. 212 (1974/75), 215--248.
 
\bibitem{Ht} R. Hartshorne, Residues and Duality, Lect. Notes Math. {
20}, Springer (1966).

\bibitem{Kon} M. Kontsevich, Homological algebra of mirror symmetry,
Proc. of the I.C.M., Vol. 1, 2 (Z{\"u}rich, 1994), 120-139,
Birkh{\"a}user, Basel,
1995, also math.AG/9411018.

\bibitem{Mu1} S. Mukai, Duality between $\D(X)$ and $\D(\hat{X})$ with its
application to Picard sheaves, Nagoya Math. J. {\bf 81} (1981) 153--175.


\bibitem{Qu} D. Quillen, Higher algebraic $K$-theory. I., in
Proc. Conf., Battelle Memorial Inst., Seattle, Wash., 1972,
Lecture Notes in Math., Vol. 341, Springer, Berlin 1973.
 




\bibitem{Ru2} A.N. Rudakov, Stability for an abelian category,
J. Algebra 197 (1997), no. 1, 231--24.

\bibitem{Sch} J.-P. Schneiders, Quasi-abelian categories and sheaves,
M{\'e}m. Soc. Math. France (NS) 76 (1999), 131p. also available at
http://www.ulg.ac.be/analg/jps/

\bibitem{Th1} R.P. Thomas, Moment maps, monodromy and mirror
manifolds,
in Symplectic Geometry and mirror symmetry: proceedings of a
conference at KIAS,
also math.AG/0212214.

\bibitem{Th2} R.P. Thomas, Stability conditions and the braid group,
Preprint math.AG/0212214.

\bibitem{Ve} J.-L.~Verdier, Des cat\'egories d\'eriv\'ees des cat\'egories
ab\'eliennes, Ast\'erisque {239} (1996).

\end{thebibliography}
\end{document}